\documentclass{amsart}                  
\usepackage{amssymb}
\usepackage{amsthm}
\newcommand{\loc}{{\rm loc}}
\newcommand{\erfc}{{\rm erfc}}
\newtheorem{theorem}{Theorem}[section]
\newtheorem{lemma}[theorem]{Lemma}
\newtheorem{corollary}[theorem]{Corollary}
\theoremstyle{remark}
\newtheorem{remark}{Remark}[section]
\begin{document}

\title[Asymptotic arbitrage and num\'eraire portfolios]{Asymptotic arbitrage and num\'eraire portfolios in large financial markets}

\author{Dmitry B. Rokhlin}


\address{D.B. Rokhlin,
Faculty of Mathematics, Mechanics and Computer Sciences,
              Southern Federal University, 
Mil'chakova str., 8a, 344090, Rostov-on-Don, Russia}
\email{rokhlin@math.rsu.ru}                   


\begin{abstract} This paper deals with the notion of a large financial market and the concepts of asymptotic arbitrage and strong asymptotic arbitrage (both of the first kind), introduced in \cite{KabKra94}, \cite{KabKra98}.
We show that the arbitrage properties of a large market are completely determined by the asymptotic behavior of the sequence of the num\'eraire portfolios, related to small markets.
The obtained criteria can be expressed in terms of contiguity, entire separation and Hellinger integrals, provided these notions are extended to sub-probability measures. 
As examples we consider market models on finite probability spaces, semimartingale and diffusion models. Also a 
discrete-time infinite horizon market model with one log-normal stock is examined. 
\end{abstract}

\maketitle

\section{Introduction} \label{intro}
The notion of a large financial market as a sequence of the traditional market models with a finite number of risky assets (called \emph{stocks} in the sequel), was introduced in \cite{KabKra94}. 
It is assumed that in $n$-th small market the discounted stock prices are described by a vector semimartingale $S^n_t=(S^{1,n}_t,\dots,S^{d(n),n}_t)$, $t\in [0,T(n)]$. Any element of the set $\mathcal X^n$ of nonnegative value processes, generated by trading strategies, is as a sum of an initial non-random endowment and a stochastic integral with respect to $S^n$. The number $d(n)$ of stocks as well as the planning horizon $T(n)$ can increase to infinity as $n\to\infty$.

The notions of asymptotic arbitrage, introduced in \cite{KabKra94}, connected the results of modern arbitrage theory and the conclusions, obtained in the framework of Capital Asset Pricing Model (Sharpe, Lintner) and Arbitrage Pricing
Model (Ross, Huberman): see \cite{Shi99}, \cite{Kab01}. 
The present paper is aimed at further study of the conditions of {\it asymptotic arbitrage of the first kind} \cite{KabKra94} and {\it strong asymptotic arbitrage of the first kind} \cite{KabKra98}.
In what follows, the term ''first kind'' is omitted since we do not consider asymptotic arbitrage of the second kind \cite{KabKra94}. Another notions of arbitrage on a large market were introduced in \cite{Kle00},  \cite{Ras03}, \cite{Kle06}.

Following \cite{KleSch96a}, asymptotic arbitrage (resp., strong asymptotic arbitrage) can be interpreted as an opportunity of getting infinitely rich with positive probability (resp., with probability 1) by risking vanishing amount of money. 
Under the assumptions of no-arbitrage and completeness of the small markets it was proved in \cite{KabKra94} that the condition of no asymptotic arbitrage (NAA) is equivalent to the contiguity of the sequence $(\mathsf P^n)$ of original probability measures with respect to the unique sequence $(\mathsf Q^n)$ of equivalent local martingale measures. 
In the paper \cite{KleSch96a} it was shown that if the small markets are incomplete then NAA condition is equivalent to the existence of some sequence $(\mathsf Q^n)$ with the same properties. Another proofs of this result were given in \cite{KleSch96b}, \cite{KabKra98}. 
It was proved in \cite{KabKra98} that the strong asymptotic arbitrage (SAA) is equivalent to the entire separation of the sequence $(\mathsf P^n)$ and any sequence $(\mathsf Q^n)$ of equivalent local martingale measures.

Let us remember that, with some abuse of terminology, a value process $0<V^n\in\mathcal X^n$ is called the \emph{num\'eraire portfolio} if any process $X^n\in\mathcal X^n$, being expressed in the units of $V^n$, becomes a supermartingale. This notion in slightly narrow meaning (martingales instead of supermartingales) was introduced in \cite{Lon90}. We refer to \cite{Bec01}, \cite{PlaHea06}, \cite{Kar06}, \cite{ChrLar07}, \cite{Chr05} for the existence 
theorems and the properties of num\'eraire portfolios as well as for further references.
The main theme of the present paper concerns the characterization of NAA and SAA conditions in terms of the sequence $(V^n)$.

The num\'eraire portfolio possesses a number of optimality properties. Particulary, under some technical conditions, it maximizes the expected logarithmic utility and the correspondent supermartingale measure (or density process) minimizes the reverse relative entropy \cite{Bec01}. So, the process $V^n$ can be looked for as the solution of the correspondent optimization problem. In any case, $V^n$ is uniquely defined \cite{Bec01}. Thus the criteria proposed below are not purely existence results, but also computational tools for checking asymptotic arbitrage conditions.  

Note that in the case of infinite time horizon a market with finite number of stocks can be regarded as ''large'' by representing it as a sequence of ''small'' markets with finite horizons $T(n)\uparrow\infty$ and the same stocks. For such a market NAA condition is tantamount to the condition of no unbounded profit with bounded risk (NUPBR) \cite{Kar06}. It was established in \cite{Kar06} that NUPBR condition is equivalent to the existence of the non-exploding num\'eraire portfolio $V$, i.e. $V_\infty<+\infty$. Underlying the connection with the results of the present paper, one may say that the condition $V_\infty<+\infty$ is imposed on the sequence $(V^n)$ of the num\'eraire portfolios, which are the restrictions of $V$ to $[0,T(n)]$. 

Somewhat surprisingly, the existence of the num\'eraire portfolios appears to be the only non-trivial assumption, concerning the structure of the small markets, allowing for the mentioned characterization of NAA and SAA conditions (see Sect.~\ref{sec:2}).
At the same time this assumption is not restrictive since in the traditional semimartingale market model with finite number of stocks and finite time horizon the existence of the num\'eraire portfolio is implied by the existence of an equivalent local martingale (or even $\sigma$-martingale) measure for the price process $(S^{1,n},\dots,S^{d,n})$ (see \cite{ChrLar07}, \cite{Kar06}).

\paragraph{Organization of the paper} In Sect.~\ref{sec:2} under minimal assumptions regarding the structure of the sets $\mathcal X^n$,
we prove that the realization of NAA and SAA conditions is completely determined by the behavior of the sequence $(V^n)$. The correspondent criteria can be expressed in terms of contiguity, entire separation and Hellinger integrals, provided these notions are extended to sub-probability measures (Theorems \ref{theo:1}, \ref{theo:2}, \ref{theo:3}, \ref{theo:4}). 

The rest of the paper is devoted to more concrete market models. Although the results of Sect.~\ref{sec:2} are applicable to all of them, the specific features of these models deserve a separate study.
In Sect.~\ref{sec:3} we consider a sequence of incomplete markets on finite probability spaces. Theorem \ref{theo:5} contains the assertion about the equivalence of NAA condition to the contiguity of the sequence $(\mathsf P^n)$ with respect to the sequence $(\widehat{\mathsf Q}^n)$ of martingale measures, minimizing the reverse relative entropy. Likewise, SAA condition is equivalent to the entire separation of $(\mathsf P^n)$ and $(\widehat{\mathsf Q}^n)$.

The results of \cite{KleSch96a}, \cite{KleSch96b}, \cite{KabKra98}, cited above and concerning the semimartingale market models, are reproved in Sect.~\ref{sec:4}. Although these results are not contained in theorems of Sect.~\ref{sec:2}, it is straightforward to give their proofs utilizing the mentioned theorems and the assertions of \cite{KraSch99}, \cite{DelSch98}, \cite{KabStr01} regarding the structure of the set of supermartingale measures for $\mathcal X^n$. This is done in Theorem \ref{theo:6}. Criteria, linking NAA and SAA conditions to the sequence of optimization problems are given in Theorem \ref{theo:7}.

In Sect.~\ref{sec:5} we treat the diffusion market models. The obtained criteria (Theorem \ref{theo:8}) are, in fact, of the same form as in \cite{KabKra98}. However, the class of models under consideration is wider since we do not require the existence of local martingale measures in the small markets.

In Sect.~\ref{sec:6} we consider a market model with discrete time and infinite horizon. It is assumed that there is only one stock with independent log-normal returns. In this example exactly one of the conditions NAA or SAA is realized. 
Under an additional assumption we are able to express these conditions in terms of convergence of some series, determined by the parameters of the model (Theorem \ref{theo:9}).

Let's briefly mention the mathematical tools used in the paper. The argumentation of Sect.~\ref{sec:2}, where the key results are collected, is based only on elementary probabilistic inequalities. In the subsequent sections
we utilize non-trivial but well-known results related to the theory of ''small'' markets. Some material from stochastic analysis is used in Sect.~\ref{sec:4} and ~\ref{sec:5}.

\section{Main results}
\label{sec:2}
Consider a sequence of probability spaces $(\Omega^n,\mathcal F^n_T,\mathsf P^n)_{n=1}^\infty$,
endowed with the filtrations $\mathbb F^n=(\mathcal F_t^n)_{t\in\mathbb T^n}$, where $\mathbb T^n$ is an interval  $[0,T(n)]$ or a set of integers $\{0,\dots,T(n)\}$. Assume that the $\sigma$-algebra $\mathcal F_0^n$ is trivial up to $\mathsf P^n$-null sets. Denote by $\mathcal X^n$ a family  of \emph{non-negative} $\mathbb F^n$-adapted stochastic processes satisfying the following conditions:
\begin{itemize}
\item[(i)] $1\in\mathcal X^n$ and $\mathcal X^n$ is a cone: if $X\in\mathcal X^n$
and $\lambda>0$ then $\lambda X\in\mathcal X^n$;
\item[(ii)] there exists a strictly positive process (\emph{num\'eraire portfolio}) 
$$V^n\in\mathcal X^n_1=\{X^n\in\mathcal X^n:X^n_0=1\}$$ such that $X^n/V^n$ is a 
$\mathsf P^n$-supermartingale  for all $X^n\in\mathcal X^n$.
\end{itemize}
The set $\mathcal X^n$ describes the value processes, generated by investment strategies in $n$-th small market. 
The large market under consideration is the sequence $(\mathcal X^n)_{n=1}^\infty$ of small markets.

Following \cite{KabKra94}, \cite{KabKra98}, we say that 
\begin{itemize}
\item
there is \emph{no asymptotic arbitrage} (NAA) on the large market $(\mathcal X^n)$ if for any sequence $X^n\in\mathcal X^n$ the condition $X^n_0\to 0$ implies that 
$$ \limsup\limits_{n\to\infty} \mathsf P^n(X^n_T\ge 1)=0;$$
\item
there exists a \emph{strong asymptotic arbitrage} (SAA) on the large market $(\mathcal X^n)$ if
$$\limsup\limits_{n\to\infty} \mathsf P^n(X^n_T\ge 1)=1$$
for some sequence $X^n\in\mathcal X^n$ such that $X^n_0\to 0$. 
\end{itemize}
By an appropriate scaling of the sequences $X^n$, it is possible to replace the sets $\{X^n_T\ge 1\}$ by the sets $\{X^n_T\ge a_n\}$, $a_n\to+\infty$, thus fitting the above definitions to the interpretation of asymptotic arbitrage, mentioned in Sect. \ref{intro}. 

In the papers \cite{KabKra94}, \cite{KleSch96a}, \cite{KleSch96b}, \cite{KabKra98} criteria for NAA and SAA conditions are expressed in terms of  contiguity and entire separation of some sequences of probability measures. In the context of the present paper it is convenient to extend these notions (see \cite{JacShi03}) to \emph{sub-probability measures}, that is to countably additive measures $\mathsf Q_n$, satisfying the condition $0\le \mathsf Q_n(A_n)\le 1$, $A_n\in\mathcal F^n_T$.

A sequence of probability measures $\mathsf P^n$ is called \emph{contiguous} with respect to a sequence of sub-probability measures $\mathsf Q^n$ (notation: $(\mathsf P^n)\lhd (\mathsf Q^n)$) if the condition $\mathsf Q^n(A^n)\to 0$, $A^n\in\mathcal F^n_T$ implies that $\mathsf P^n(A^n)\to 0$.
A sequence of probability measures $\mathsf P^n$ and a sequence of sub-probability measures $\mathsf Q^n$ are called
\emph{entirely (asymptotically) separated} (notation: $(\mathsf P^n)\bigtriangleup (\mathsf Q^n)$) if there exist a sequence of natural numbers $n_k\uparrow\infty$ and sets $A^{n_k}\in\mathcal F^{n_k}_T$ such that $\mathsf P^{n_k}(A^{n_k})\to 1$ and $\mathsf Q^{n_k}(A^{n_k})\to 0$, $k\to\infty$.

We say that a strictly positive $\mathbb F^n$-adapted stochastic process $Z^n$ is an {\it equivalent supermartingale density} for $\mathcal X^n$ if $X^nZ^n$ is a $\mathsf P^n$-supermartingale for all
$X^n\in\mathcal X^n$ and $Z^n_0=1$. Denote by $\mathcal D^n$ the set of all equivalent supermartingale densities (compare with \cite{KraSch99}, \cite{Zit02},  \cite{Kar06}). Clearly, $1/V^n\in\mathcal D^n$.

Denote by $L^{0,n}_+=L^0(\Omega^n,\mathcal F^n_T,\mathsf P^n)$ the set of (equivalence classes) of non-negative $\mathcal F^n_T$-measurable random variables. Let $\xi^n\in L^{0,n}_+$. Following \cite{Kal97}, we denote by $\xi^n\cdot\mathsf P^n$ the measure with the $\mathsf P^n$-density $\xi^n$: 
$$(\xi^n\cdot\mathsf P^n)(A^n)=\mathsf E_{\mathsf P^n}(\xi^n I_{A^n}),\ \ \ A^n\in\mathcal F^n_T.$$ 
To each process $Z^n\in\mathcal D^n$ we assign a sub-probability measure $Z_T^n\cdot\mathsf P^n$. 

Consider a sequence $\xi^n\in L^{0,n}_+$ and a sequence of probability measures $\mathsf P^n$.  
According to the definition of \cite{JacShi03} a sequence  $(\xi^n|\mathsf P^n)$ \emph{is tight} if
$$	\lim_{M\to\infty}\limsup_{n\to\infty}\mathsf P^n(\xi^n\ge M)=0. $$
\begin{theorem} \label{theo:1}
The following conditions are equivalent:
\begin{itemize}
\item[(a)] NAA;
\item[(b)] $(Y^n_T|\mathsf P^n)$ is tight for any sequence $Y^n\in\mathcal X^n_1$;
\item[(c)] $(V^n_T|\mathsf P^n)$ is tight;
\item[(d)] $(\mathsf P^n)\lhd ((V_T^n)^{-1}\cdot\mathsf P^n)$;
\item[(e)] $(\mathsf P^n)\lhd (Z_T^n\cdot\mathsf P^n)$ for some sequence $Z^n\in\mathcal D^n$.
\end{itemize}
\end{theorem}
\textit{Proof.} (a) $\Longrightarrow$ (b). Assume that condition (b) is violated.
Then there exist a number $\beta>0$ and a sequence $Y^n\in\mathcal X^n_1$ such that
$$ \limsup_{n\to\infty}
    \mathsf P^n(Y^n_T\ge M)\ge\beta>0$$
for all $M>0$. Take an increasing sequence of natural numbers $n_k$ so as
$\mathsf P^{n_k}(Y^{n_k}_T\ge k)\ge\beta/2$.
Further, take a sequence $M^n\uparrow\infty$ meeting the condition $M^{n_k}=k$ and put $X^n=Y^n/M^n$. Then $X^n_0=1/M^n\to 0$ and
$$ \limsup_{n\to\infty} \mathsf P^n(X^n_T\ge 1)\ge
   \limsup_{k\to\infty} \mathsf P^{n_k}(Y^{n_k}_T\ge k)\ge
   \beta/2>0.$$
This means that NAA condition is violated. 

(c) $\Longrightarrow$ (d). Let $((V_T^n)^{-1}\cdot\mathsf P^n)(A^n)\to 0$,
$A^n\in\mathcal F^n_T$. Then the inequality
\begin{eqnarray*}
\mathsf P^n(A^n) &= & \mathsf P^n(A^n\cap\{V^n_T\ge M\})+
\int\limits_{A^n\cap\{V^n_T< M\}} V^n_T d((V_T^n)^{-1}\cdot\mathsf P^n)\\
&\le & \mathsf P^n(V^n_T\ge M)+M((V_T^n)^{-1}\cdot\mathsf P^n)(A^n)
\end{eqnarray*}
yields that 
$$ \limsup_{n\to\infty} \mathsf P^n(A^n)\le
   \limsup_{n\to\infty} \mathsf P^n(V^n_T\ge M).$$
By condition (c) it follows that $\mathsf P^n(A^n)\to 0$.

(e) $\Longrightarrow$ (a). Consider a sequence $X^n\in\mathcal X^n$ such that $\lim_{n\to\infty} X^n_0\to 0$. 
The inequality
\begin{equation} \label{1.1}
(Z^n_T\cdot\mathsf P^n)(X^n_T\ge 1)=\mathsf E_{\mathsf P^n}\left(Z^n_T X^n_T 
\frac{I_{\{X^n_T\ge 1\}}}{X^n_T}\right)\le\mathsf E_{\mathsf P^n}(Z^n_T X^n_T )\le X^n_0
\end{equation}
and the contiguity of $(\mathsf P^n)$ with respect to $(Z^n_T\cdot\mathsf P^n)$ imply that 
$$\lim_{n\to\infty}\mathsf P^n(X^n_T\ge 1)=0,$$
i.e. NAA condition is satisfied.

Implications (b)  $\Longrightarrow$ (c) and (d) $\Longrightarrow$ (e) are evident. \qed

\begin{theorem} \label{theo:2}
The following conditions are equivalent:
\begin{itemize}
\item[(a)] SAA;
\item[(b)] there exists a sequence $Y^n\in\mathcal X^n_1$ such that
$$ \limsup_{n\to\infty}\mathsf P^n(Y^n_T\ge M)=1\ \textnormal{for all}\ M>0;$$
\item[(c)] $\limsup_{n\to\infty}\mathsf P^n(V^n_T\ge M)=1$ for all $M>0$;
\item[(d)] $(\mathsf P^n)\bigtriangleup ((V_T^n)^{-1}\cdot\mathsf P^n)$;
\item[(e)] $(\mathsf P^n)\bigtriangleup (Z_T^n\cdot\mathsf P^n)$ for any sequence $Z^n\in\mathcal D^n$.
\end{itemize}
\end{theorem}
\textit{Proof.} 
(b) $\Longrightarrow$ (a).
Take an increasing sequence of natural numbers $n_k$ such that
$$\mathsf P^{n_k}(Y^{n_k}_T\ge k)\ge 1-1/k.$$
Define the sequences $M_n$, $X^n$ as in the proof of the implication (a) $\Longrightarrow$ (b) of Theorem~\ref{theo:1}. 
We have $X^n_0\to 0$ and 
$$\limsup_{n\to\infty}\mathsf P^n(X^n_T\ge 1)\ge
  \limsup_{k\to\infty}\mathsf P^{n_k}(X^{n_k}_T\ge 1)=
  \limsup_{k\to\infty}\mathsf P^{n_k}(Y^{n_k}_T\ge k)=1.$$
Thus, $(X^n)$ realizes the strong asymptotic arbitrage.

(d) $\Longrightarrow$ (c). Let the sequences 
$n_k\uparrow\infty$ and $A^{n_k}\in\mathcal F^{n_k}_T$ be such that
$$\mathsf P^{n_k}(A^{n_k})\to 1,\ \ \ ((V_T^{n_k})^{-1}\cdot\mathsf P^{n_k})(A^{n_k})\to 0.$$
Then
\begin{eqnarray*}
\mathsf P^{n_k}(V^{n_k}_T<M) &=&
\mathsf P^{n_k}(\{V^{n_k}_T<M\}\cap (\Omega\backslash A^{n_k}))\\
&+& \int\limits_{\{V^{n_k}_T<M\}\cap A^{n_k}} V^{n_k}_T d((V_T^{n_k})^{-1}\cdot\mathsf P^{n_k}) \\
&\le & \mathsf P^{n_k}(\Omega\backslash A^{n_k})
+ M ((V_T^{n_k})^{-1}\cdot\mathsf P^{n_k})(A^{n_k})\to 0,\ \ k\to\infty
\end{eqnarray*}
and condition (c) is satisfied.

(a) $\Longrightarrow$ (e). Let $X^n$ be a sequence, mentioned in the definition of the strong asymptotic arbitrage and $Z^n\in\mathcal D^n$. The inequality (\ref{1.1}) and the definition of the strong asymptotic arbitrage show that 
$$ \lim_{n\to\infty} (Z_T^n\cdot\mathsf P^n)(X^n_T\ge 1)=0,\ \ \ 
   \limsup_{n\to\infty} \mathsf P^n(X^n_T\ge 1)=1.$$
Therefore the sequences $(\mathsf P^n)$ and $(Z_T^n\cdot\mathsf P^n)$ are entirely separated.

Implications (c) $\Longrightarrow$ (b) and (e) $\Longrightarrow$ (d) are evident. \qed

Inspired by the results of \cite{JacShi03}, \cite{KabKra98}, we derive criteria for NAA and SAA conditions in terms of Hellinger-type integrals. Auxilary inequalities, suitable for this purpose, are collected in Lemma \ref{lem:1}. Inequality (\ref{1.4}) is applied also in Sect.~\ref{sec:5}. The proof of inequality (\ref{1.3}) is, in fact, borrowed from \cite{JacShi03} (p. 287). 
\begin{lemma} \label{lem:1}
Let $\xi$ be a non-negative random variable, defined on a probability space $(\Omega,\mathcal F,\mathsf P)$ and
let $\alpha\in (0,1)$, $M>0$, $N>0$ be some numbers. Then
\begin{equation} \label{1.2}
	\mathsf P(\xi<M)\le M^\alpha\mathsf E \xi^{-\alpha};
\end{equation}
if $\mathsf E(1/\xi)\le 1$ then
\begin{equation} \label{1.3}
	\mathsf E \xi^{-\alpha}\le \left(\frac{1}{M}\right)^\alpha+\left(\frac{1}{N}\right)^{1-\alpha}+
   N^\alpha\mathsf P(\xi<M);
\end{equation}
if a non-negative random variable $\eta$ is such that $\mathsf E(\xi/\eta)\le 1$, then
\begin{equation}  \label{1.4}
	\mathsf P(\xi\ge M)\le\frac{N}{M}+\mathsf P(\eta\ge N).
\end{equation}
\end{lemma}
\textit{Proof.} The estimate (\ref{1.2}) follows from the inequality
$$I_{\{\xi<M\}}\le M^\alpha \xi^{-\alpha} I_{\{\xi<M\}}.$$

Now let $\mathsf E (1/\xi)\le 1$. Then
\begin{eqnarray*}
 && \mathsf E \xi^{-\alpha} =
   \mathsf E (\xi^{-\alpha}I_{\{\xi\ge M\}})+
   \mathsf E \left(\frac{\xi^{1-\alpha}}{\xi}
    I_{\{\xi< M\}}I_{\{\xi<1/N\}}\right)\nonumber\\
   &+&\mathsf E(\xi^{-\alpha}I_{\{\xi<M\}}
   I_{\{\xi\ge 1/N\}}) \le \left(\frac{1}{M}\right)^\alpha+\left(\frac{1}{N}\right)^{1-\alpha}+
   N^\alpha\mathsf P(\xi<M).
\end{eqnarray*}

At last, assume that non-negative random variables $\xi$, $\eta$ verify the inequality  $\mathsf E(\xi/\eta)\le 1$. We see that
\begin{equation}\label{1.5}
\mathsf E\left(\frac{I_{\{\xi\ge M\}}}{\eta}\right)\le\mathsf E\left(\frac{\xi}{M}
\frac{I_{\{\xi\ge  M\}}}{\eta}\right)\le\frac{1}{M}.
\end{equation}
\begin{equation}\label{1.6}
\mathsf E\left(\frac{I_{\{\xi\ge M\}}}{\eta}\right)\ge
\mathsf E\left(\frac{I_{\{\xi\ge M\}} I_{\{\eta<N\}}}{\eta}\right)
\ge\frac{1}{N}\mathsf P(\{\xi\ge M\}\cap\{\eta<N\})
\end{equation}
But for any sets $B_1, B_2\in \mathcal F$ we have
$$ \mathsf P(B_1\cap B_2)=\mathsf P(B_1)+\mathsf P(B_2)-
   \mathsf P(B_1\cup B_2)\ge \mathsf P(B_1)-
   \mathsf P(\Omega\backslash B_2).$$
Hence (\ref{1.5}) and (\ref{1.6}) imply (\ref{1.4}):
$$\frac{N}{M}\ge \mathsf P(\{\xi\ge M\}\cap\{\eta<N\})\ge
\mathsf P(\xi\ge M)-\mathsf P(\eta\ge N). \ \ \ \ \qed $$

\begin{corollary} \label{cor:1}
Any sequence $Z^n\in\mathcal D^n$ satisfies the equality
\begin{equation} \label{1.7}
	\lim_{\alpha\downarrow 0}\liminf_{n\to\infty}\mathsf E_{\mathsf P_n} (Z^n_T)^\alpha
=\lim_{M\to\infty}\liminf_{n\to\infty}\mathsf P^n((Z^n_T)^{-1}<M).
\end{equation}
In particular,  $((Z^n_T)^{-1}|\mathsf P^n)$ is tight iff the left-hand side of (\ref{1.7}) is equal to 1.
\end{corollary}
\textit{Proof.} Applying inequality (\ref{1.2}):
$$\mathsf P^n((Z_T^n)^{-1}<M)\le M^\alpha\mathsf E_{\mathsf P^n} (Z^n_T)^\alpha $$
and taking the limits as $n\to\infty$, $\alpha\to 0$ and $M\to\infty$, we get
$$ \lim_{M\to\infty}\liminf_{n\to\infty}\mathsf P^n((Z^n_T)^{-1}<M)\le
   \lim_{\alpha\downarrow 0}\liminf_{n\to\infty}\mathsf E_{\mathsf P_n} (Z^n_T)^\alpha.$$
Furthermore, using (\ref{1.3}):
$$ \mathsf E_{\mathsf P^n} (Z^n_T)^\alpha\le \left(\frac{1}{M}\right)^\alpha+\left(\frac{1}{N}\right)^{1-\alpha}+
   N^\alpha\mathsf P^n((Z^n_T)^{-1}<M),$$
and taking the limits as $n\to\infty$, $M\to\infty$, $\alpha\to 0$ and $N\to\infty$, we obtain the reverse inequality. \qed

\begin{corollary} \label{cor:2}
The following conditions are equivalent:
\begin{itemize}
  \item[(a)] $\limsup\limits_{n\to\infty}\mathsf P^n(V^n_T\ge M)=1$ for all $M>0$;
  \item[(b)] $\liminf\limits_{n\to\infty}\mathsf E_{\mathsf P_n} (V^n_T)^{-\alpha}=0$ for some $\alpha\in (0,1)$.
\end{itemize}
\end{corollary} 
\textit{Proof.} Implication (b) $\Longrightarrow$ (a) follows from (\ref{1.2}):
$$ \mathsf P(V^n_T<M)\le M^\alpha\mathsf E (V^n_T)^{-\alpha}.$$
Conversely, if (a) holds then by (\ref{1.3}) we get
$$ \liminf_{n\to\infty}	\mathsf E (V^n_T)^{-\alpha}\le \left(\frac{1}{M}\right)^\alpha+\left(\frac{1}{N}\right)^{1-\alpha}.$$
Since $M$, $N$ are arbitrary, this yields (b). \qed

\begin{corollary} \label{cor:3}
Let $X^n\in\mathcal X^n_1$, $Z^n\in\mathcal D^n$. Then
$$\limsup_{n\to\infty}\mathsf P^n(X^n_T\ge M)\le\frac{N}{M}+\limsup_{n\to\infty}\mathsf P^n(V^n_T\ge N);$$
$$\limsup_{n\to\infty}\mathsf P^n(V^n_T\ge M)\le\frac{N}{M}+\limsup_{n\to\infty}\mathsf P^n(1/Z^n_T\ge N).$$
In particular, if $((Z^n_T)^{-1}|\mathsf P^n)$ is tight then $(V^n_T|\mathsf P^n)$ is tight as well.
\end{corollary} 
The proof directly follows from (\ref{1.4}) since 
$$\mathsf E_{\mathsf P^n}(X^n_T/V^n_T)\le 1,\ \ \ \mathsf E_{\mathsf P^n}(V^n_T/(Z^n_T)^{-1})\le 1. \ \ \ \qed$$

\begin{theorem} \label{theo:3}
The following conditions are equivalent:
\begin{itemize}
\item[(a)] NAA;
\item[(b)] $\lim\limits_{\alpha\downarrow 0}\liminf\limits_{n\to\infty}
\mathsf E_{\mathsf P_n} (V^n_T)^{-\alpha}=1$;
\item[(c)] there exists a sequence $Z^n\in\mathcal D^n$ such that
$$\lim_{\alpha\downarrow 0}\liminf_{n\to\infty}\mathsf E_{\mathsf P_n} (Z^n_T)^\alpha=1.$$ 
\end{itemize}
\end{theorem}
\textit{Proof.} (a) $\iff$ (b). By Theorem \ref{theo:1} the absence of asymptotic arbitrage is equivalent to the tightness of the sequence $(V^n_T|\mathsf P^n)$. But in accordance with Corollary \ref{cor:1} this condition is equivalent to (b) since $(V^n)^{-1}\in\mathcal D^n$. 

(c) $\Longrightarrow$ (a). Corollaries \ref{cor:1} and \ref{cor:3} imply that $(V^n_T|\mathsf P^n)$ is tight.
Applying again Theorem~\ref{theo:1}, we conclude that there is no asymptotic arbitrage. 
The last implication (b) $\Longrightarrow$ (c) is evident. \qed

\begin{theorem} \label{theo:4}
The following conditions are equivalent:
\begin{itemize}
\item[(a)] SAA;
\item[(b)] $\liminf\limits_{n\to\infty}\mathsf E_{\mathsf P_n} (V^n_T)^{-\alpha}=0$ for some $\alpha\in (0,1)$;
\item[(c)] the exists a sequence $Y^n\in\mathcal X^n_1$ such that
$$\liminf_{n\to\infty}\mathsf E_{\mathsf P_n} (Y^n_T)^{-\alpha}=0\ \textit{for\ some}\ \alpha\in (0,1).$$
\end{itemize}
\end{theorem}
\textit{Proof.} The assertion about the equivalence of (a) and (b) follows from Theorem~\ref{theo:2} and 
Corollary~\ref{cor:2}.

(c) $\Longrightarrow$ (a). Under condition (c) inequality (\ref{1.2}) yields:
$$ \liminf_{n\to\infty}\mathsf P^n(Y^n_T <M)=0.$$
Therefore condition (b) of Theorem~\ref{theo:2} holds true.
Implication (b) $\Longrightarrow$ (c) is evident. \qed

\begin{remark} \label{rem:1} 
It is equivalent to require that conditions (b) of Corollary~\ref{cor:2} and Theorem~\ref{theo:4} hold true for all $\alpha\in (0,1)$ or in the limit as $\alpha\downarrow 0$ (compare with \cite{JacShi03}, p.287).
\end{remark}

\begin{remark} \label{rem:2}
Extending the well-known terminology \cite{JacShi03} to sub-probability measures, the expression
$$ H(\alpha;Z^n_T\cdot\mathsf P^n,\mathsf P^n)=
  \mathsf E_{\mathsf P^n}\left(\frac{d(Z^n_T\cdot\mathsf P^n)}{d\mathsf P^n}\right)^\alpha=
  \mathsf E_{\mathsf P_n}(Z^n_T)^\alpha$$
can be called a \emph{Hellinger integral} of order $\alpha\in (0,1)$ between $Z^n_T\cdot\mathsf P^n$ and $\mathsf P^n$. The related interpretation of Theorems \ref{theo:3}, \ref{theo:4} can be compared with \cite{KabKra98}. 
\end{remark}

\begin{remark} \label{rem:3} It follows from Jensen's inequality and the definition of $V^n$ that the latter process is \emph{relatively log-optimal} \cite{ChrLar07}, \cite{Kar06}:
$$\mathsf E_{\mathsf P^n}\left(\ln\frac{X^n_T}{V^n_T}\right)\le 0,\ \ \ X^n\in\mathcal X^n_1.$$
If, moreover, $\mathsf E_{\mathsf P^n}(\ln V^n_T)<\infty$ then $V^n$ is {\it log-optimal} and the supermartingale density $\widehat Z^n=1/V^n$ minimizes the \emph{reverse relative entropy} among all supermartingale
densities \cite{Bec01}, \cite{Kar06}. In other words,
$$\sup_{X\in\mathcal X^n_1}\mathsf E_{\mathsf P^n}\ln X^n_T=\mathsf E_{\mathsf P^n}(\ln V^n_T)=
  \mathsf E_{\mathsf P^n}\left(\ln\frac{1}{\widehat Z_T^n}\right)= 
  \inf_{Z\in\mathcal D^n}\mathsf E_{\mathsf P^n}\left(\ln\frac{1}{Z_T^n}\right).$$
It should be mentioned that if $Z^n\cdot\mathsf P^n$ is a probability measure then $\mathsf E_{\mathsf P^n}\ln(1/Z^n_T)$ coincides with the entropy of $\mathsf P^n$ with respect to $Z^n\cdot\mathsf P^n$ in the usual meaning (see e.g. \cite{CovTho06}).
\end{remark}

\section{Market models on finite probability spaces}
\label{sec:3}
Assume that for every $n$ the set $\Omega^n$ is finite, $\mathbb T^n=\{0,\dots,T(n)\}$ and there is an $\mathbb F^n$-adapted
process $S^n_t=(S^{1,n}_t,\dots,S^{d(n),n}_t)$, $t\in\mathbb T^n$ of discounted stock prices, defined on the probability space $(\Omega^n,\mathcal F^n_T,\mathsf P^n)$. 
Denote by $(x,y)$ the scalar product of vectors $x$ and $y$ in $\mathbb R^d$.
Suppose that the set $\mathcal X^n$ of value processes consists of the elements $X^n$, admitting the representation
$$ X_t^n=x^n+\sum_{j=1}^t(\gamma_j^n, S_j^n-S_{j-1}^n)\ge 0,\ \  
t=0,\dots, T(n),$$
where $x^n\in\mathbb R$ and the components of $\mathbb F^n$-predictable processes $$\gamma^n=(\gamma^{1,n}, \dots, \gamma^{d(n),n})$$ describe the number of stocks in investor's portfolio.

Let $\mathbf M^n$ be the set of equivalent to $\mathsf P^n$ probability measures, under which the
process $S^n$ is a martingale. 
Denote by $H(\mathsf P|\mathsf Q)=\mathsf E_\mathsf P \ln\left(d\mathsf P/d\mathsf Q\right)$ 
the entropy of $\mathsf P$ with respect to an (equivalent) measure $\mathsf Q$.

\begin{theorem} \label{theo:5}
 Let $\mathbf M^n\neq\emptyset$. Then for every $n$ there exists a unique minimal reverse entropy martingale measure
$\widehat{\mathsf Q}^n\in\mathbf M^n$:
$$ H(\mathsf P^n|\widehat{\mathsf Q}^n)\le H(\mathsf P^n|\mathsf Q^n),\ \ \ \mathsf Q^n\in\mathbf M^n$$
and the following relations hold true
$$ NAA \iff (\mathsf P^n)\lhd (\widehat{\mathsf Q}^n);\ \ \ SAA \iff (\mathsf P^n)\bigtriangleup (\widehat{\mathsf Q}^n). $$
\end{theorem}
\textit{Proof.} Consider the optimization problems
\begin{equation} \label{3.1}
\mathsf E_{\mathsf P^n}\ln X^n_T\mapsto\max_{X^n\in\mathcal X^n_1},\ \ \ \mathsf E_{\mathsf P^n}\ln\frac{d\mathsf P^n}{d\mathsf Q^n}\mapsto\min_{\mathsf Q^n\in\mathbf M^n}.
\end{equation}
From Theorem 2.4 of \cite{Sch04} we know that the problems (\ref{3.1}) are solvable and their unique solutions $\widehat V^n\in\mathcal X^n_1$, $\widehat{\mathsf Q}^n\in\mathbf M^n$ are related as follows $\widehat{\mathsf Q}^n=(\widehat V^n_T)^{-1}\cdot \mathsf P^n$. Moreover, the process $\widehat V^n=V^n$ is the num\'eraire portfolio and $\widehat{\mathsf Q}^n$  minimizes the reverse relative entropy (see \cite{Bec01} and Remark \ref{rem:3}). Thus the proof is implied by conditions (d) of Theorems \ref{theo:1} and \ref{theo:2}. \qed

\section{Semimartingale market models}
\label{sec:4}
In this section we assume that the processes in $\mathcal X^n$, $\mathcal D^n$ have c\'adl\'ag trajectories $\mathsf P^n$-a.s. Evidently, these assumptions do not affect the argumentation of Sect.~\ref{sec:2}.

Suppose there are $d(n)$ stocks in $n$-th small market and their discounted prices are described by a vector 
semimartingale $S^n=(S^{1,n},\dots,S^{d(n),n})$, adapted to the filtration $\mathbb F^n=(\mathcal F^n_t)_{0\le t\le T(n)}$, satisfying the usual conditions \cite{JacShi03}.
Furthermore, assume that any element $X^n\in\mathcal X^n$ is of the form 
$$X_t^n=x^n+(\gamma^n\circ S^n)_t\ge 0,\ \ t\in [0,T(n)],$$ 
where $x^n\in\mathbb R$, $\gamma^n=(\gamma^{1,n},\dots,\gamma^{d(n),n})$ is an $S^n$-integrable process (notation: $\gamma^n\in L(S^n)$) and
$$(\gamma^n\circ S^n)_t=\int_{0}^t(\gamma^n_u,dS^n_u)$$
is a vector stochastic integral \cite{JacShi03}. The quantity $\gamma_t^{i,n}$ determines the number of units of $i$-th stock in investor's portfolio at time $t$.


Let $\mathbf M_\sigma^n$ (resp., $\mathbf M_\loc^n$) be the set of equivalent to $\mathsf P^n$ probability measures, under which the process $S^n$ is a $\sigma$-martingale (resp., local martingale). 
Assume that $\mathbf M_\sigma^n\neq\emptyset$. Then there exists
the num\'eraire portfolio $V^n\in\mathcal X^n_1$. In full generality this assertion follows from the results of \cite{ChrLar07} (see also \cite{Bec01}, \cite{Kar06}). Thus the conclusions of Sect.~\ref{sec:2} are valid for the sequence $(\mathcal X^n)$ of semimartingale market models.

Denote by $\mathbf M_s^n$ the set of equivalent to $\mathsf P^n$ probability measures $\mathsf Q^n$ such that all processes in $\mathcal X^n$ are $\mathsf Q^n$-supermartingales and put $\mathbf D^n=\{Z_T^n\cdot\mathsf P^n:Z^n\in\mathcal D^n\}$.

The statements of Theorem \ref{theo:6} below were proved in \cite{KabKra98} under the assumption $\mathbf M_\loc^n\neq\emptyset$. 
Previously, assertion (a) was proved in \cite{KleSch96a}, where the process $S^n$ was assumed to be locally bounded.
We give alternative proofs, based on Theorems \ref{theo:1}, \ref{theo:2} and on the well-known non-trivial results,
concerning the approximative properties of inclusions
$$\mathbf M_\loc^n\subset\mathbf M_\sigma^n\subset\mathbf M_s^n\subset\mathbf D^n.$$
Note that the inclusion $\mathbf M_\sigma^n\subset\mathbf M_s^n$ follows from the Ansel-Stricker theorem  \cite{AnsStr94} (Corollary 3.5): see \cite{DelSch98} (Theorem 5.3).

\begin{theorem} \label{theo:6}
Let $\mathbf M_\sigma^n\neq\emptyset$. Then 
\begin{itemize}
	\item[(a)]  NAA condition is satisfied iff there exists a sequence $\mathsf Q^n\in \mathbf M_\sigma^n$ such that $(\mathsf P^n)\lhd(\mathsf Q^n)$;
	\item[(b)] SAA condition is satisfied iff $(\mathsf P^n) \bigtriangleup (\mathsf Q^n)$ for any sequence $\mathsf Q^n\in \mathbf M_\sigma^n$.
\end{itemize}
If $\mathbf M_{\loc}^n\neq\emptyset$ then in the above statements the set $\mathbf M_\sigma^n$ can be replaced by
$\mathbf M_{\loc}^n$.
\end{theorem}
\textit{Proof.} Let $Z^n\in\mathcal D^n$. It follows from Theorem 4(b) of \cite{Zit02} (and, in fact, from the results of \cite{KraSch99}, \cite{BraSch99}) that there exist a sequence of non-negative $\mathcal F^n_T$-measurable random variables $g^{k,n}\le 1$ and a sequence $Z^{k,n}_T\cdot\mathsf P^n\in\mathbf M_s^n$, $k\ge 1$  such that
$$g^{k,n} Z^{k,n}_T\to Z^n_T,\ \ k\to\infty,\ \ \mathsf P^n\textnormal{-a.s.}$$

To each probability measure $\mathsf Q^n$, absolutely continuous with respect to $\mathsf P^n$, we assign its Radon-Nikodym derivative $d\mathsf Q^n/d\mathsf P^n$. Thus the set of all such measures is identified with the subset $L^{1,n}_+=L^1_+(\Omega^n,\mathcal F^n_T,\mathsf P^n)$ of non-negative $\mathsf P^n$-integrable random variables. 
It is proved in \cite{DelSch98} (Proposition 4.7) that the set $\mathbf M_\sigma^n$ is dense in 
$\mathbf M_s^n$ with respect to the norm topology of $L^{1,n}$. Therefore, we may assume that $Z^{k,n}_T\cdot\mathsf P^n\in\mathbf M_\sigma^n$.

Put $\varepsilon^n=1/n$ and choose $M(n)>0$ so that 
$$\mathsf E_{\mathsf P^n} \left(Z^n_T I_{\{Z^n_T\ge M(n)\}}\right)\le\varepsilon^n.$$
By Egorov's theorem there is a set $B^n\in\mathcal F^n_T$ such that 
$$\mathsf P^n(\Omega^n\backslash B^n)\le\varepsilon^n/M(n)$$ 
and $g^{k,n} Z^{k,n}_T$ converges to $Z^n_T$ uniformly on $B^n$ as $k\to\infty$.
Furthermore, choose $k(n)$ large enough to ensure the inequality 
$$|g^{k(n),n} Z^{k(n),n}_T- Z^n_T|\le\varepsilon^n\ \textnormal{on}\ B^n.$$
Then for any $A^n\in\mathcal F^n_T$ the following esimates hold true:
$$ \int\limits_{A^n\cap (\Omega\backslash B^n)}Z^n_T\,d\mathsf P^n\le \int\limits_{\{Z^n_T\ge M(n)\}}Z^n_T\,d\mathsf P^n
+\int\limits_{\{Z^n_T<M(n)\}\cap (\Omega\backslash B^n)}Z^n_T\,d\mathsf P^n\le 2\varepsilon^n,$$
$$ \int\limits_{A^n\cap B^n}Z^n_T\,d\mathsf P^n\le
   \varepsilon^n+\int\limits_{A^n}Z^{k(n),n}_T\,d\mathsf P^n.$$
Consequently,
\begin{equation} \label{2.1}
(Z_T^n\cdot\mathsf P^n)(A^n)\le 3\varepsilon^n+(\mathsf Q^n)(A^n),
\end{equation}
where $\mathsf Q^n= Z^{k(n),n}_T\cdot\mathsf P^n$. 

It follows from inequality (\ref{2.1}) that if the sequence $(\mathsf P^n)$ is contiguous with respect to
$(Z_T^n\cdot\mathsf P^n)$ then it is contiguous with respect to some sequence $\mathsf Q^n\in\mathbf M_\sigma^n$. Likewise, if $(\mathsf P^n)\bigtriangleup (\mathsf Q^n)$ for any sequence $\mathsf Q^n\in\mathbf M_\sigma^n$ then $(\mathsf P^n)\bigtriangleup (Z_T^n\cdot\mathsf P^n)$ for any sequence $Z^n\in\mathcal D^n$. By conditions (e) of Theorems \ref{theo:1} and \ref{theo:2} this proves that NAA condition implies the existence of a sequence $\mathsf Q^n\in\mathbf M_\sigma^n$ such that $(\mathsf P^n)\lhd(\mathsf Q^n)$, and SAA condition is satisfied if the relation $(\mathsf P^n)\bigtriangleup (\mathsf Q^n)$ holds true for any sequence $\mathsf Q^n\in\mathbf M_\sigma^n$. 

The converse statements also follow from conditions (e) of Theorems \ref{theo:1}, \ref{theo:2} and the inclusion $\mathbf M_\sigma^n\subset\mathbf D^n$.

Let's prove the last assertion of the theorem. If $\mathbf M_\loc^n\neq\emptyset$ then $\mathbf M_\loc^n$ is dense in  $\mathbf M_\sigma^n$ with respect to the norm topology of $L^{1,n}$. This result, as is mentioned in \cite{Klo06}, follows from Theorem 1.1 of the paper \cite{KabStr01}. Thus in the above argumentation the set $\mathbf M_\sigma^n$ can be replaced by
$\mathbf M_{\loc}^n$. \qed

The subsequent theorem shows, in particular, that NAA and SAA conditions can be checked by exploring the sequence of strategies that achieve the optimal expected logarithmic utility of the terminal wealth. Taking into account the results of \cite{GolKal00}, \cite{GolKal03}, this gives an opportunity to express criteria for NAA and SAA conditions in terms of the semimartingale characteristic sequences of the stock price processes $S^n$.

\begin{theorem} \label{theo:7}
Suppose, $\mathbf M^n_\sigma\neq\emptyset$ and 
$$\sup_{X\in\mathcal X^n_1}\mathsf E_{\mathsf P^n}\ln X^n_T<\infty.$$
Then for every $n$ there exists a unique log-optimal strategy $\widehat V^n$, a unique minimal reverse entropy supermartingale density $\widehat Z^n$ (see Remark \ref{rem:3}) and the following relations hold true: 
$$ NAA \iff \lim_{\alpha\downarrow 0}\liminf_{n\to\infty} \mathsf E_{\mathsf P_n} (\widehat V^n_T)^{-\alpha}=1
\iff (\mathsf P^n)\lhd (\widehat Z_T^n\cdot\mathsf P^n);$$
$$ SAA \iff \liminf_{n\to\infty} \mathsf E_{\mathsf P_n} (\widehat V^n_T)^{-\alpha}=0\ \textit{for some}\ \alpha\in (0,1)
\iff (\mathsf P^n)\bigtriangleup (\widehat Z_T^n\cdot\mathsf P^n).$$
\end{theorem}
\textit{Proof.} The statements, concerning the existence and uniqueness of the processes $\widehat V^n$, $\widehat Z^n$ with the above properties, the existence of the num\'eraire portfolio $V^n$ and the relations  $V^n=\widehat V^n$, $1/V^n=\widehat Z^n$ were proved in \cite{Bec01} with the use of the results of \cite{KraSch99}. It remains to apply conditions (b) of Theorems \ref{theo:3} and \ref{theo:4} 
and conditions (d) of Theorems \ref{theo:1} and \ref{theo:2}. \qed

\section{Diffusion market models}
\label{sec:5}
In the framework of the model considered in Sect.~\ref{sec:4}, assume that for every $n$ the stock prices are driven by $m(n)$-dimensional process $(W^{1,n}_t,\dots,W^{m(n),n}_t)$, $t\in [0, T(n)]$, composed of independent standard Wiener processes $W^{i,n}$, defined on the filtered probability space $(\Omega^n,\mathcal F^n,\mathsf P^n,\mathbb F^n)$. The stock prices $S^n$ are subject to the system of stochastic differential equations
$$ dS^{i,n}_t=S^{i,n}_t dR_t^{i,n},\ \ \ R_t^{i,n}=\int_0^t\mu_u^{i,n}du+\int_0^t (\beta^{i,n}_u,d W^n_u),\ \ i=1,\dots,d(n),$$
where predictable stochastic processes $\mu^{i,n}$, $\beta^{i,n}$ satisfy the conditions
$$\int_0^{T(n)}\left(|\mu^{i,n}_t|+|\beta^{i,n}_t|^2\right)dt<+\infty$$ and $S_0^i>0$. We use the notation $|x|=\sqrt{(x,x)}$ for the length of a vector $x$ as well as for the absolute value of a scalar.

Assume also that $d(n)\le m(n)$ and the rank of the matrix $\sigma^n$ with the rows $(\beta^{i,n})_{i=1}^{d}$ is equal to $d(n)$ for all $t,\omega$ and $n$. Then the matrix $\sigma^n(\sigma^n)^T$ is invertible. Let the vector $\lambda^n=(\sigma^n)^T(\sigma^n(\sigma^n)^T)^{-1}\mu^n$ be such that 
\begin{equation}\label{4.1}
\int_0^{T(n)}|\lambda^n_t|^2\, dt<+\infty\ \ \textnormal{a.s.}
\end{equation}
It is customary to call $\lambda^n$ a market price of risk process.

The following assertion is known: see e.g. Example 4.2 in \cite{GolKal00} and Example 2.2.19 in \cite{Chr05}. For  convenience, we provide its direct proof.
\begin{lemma} \label{lem:2}
Under the adopted assumptions there exists the num\'eraire portfolio 
$$V^n_t= \exp\left(\frac{1}{2}\int_0^t |\lambda_u^n|^2\,du+\int_0^t(\lambda_u^n,dW_u^n)\right),\ \ V\in\mathcal X^n_1.$$
\end{lemma}
\emph{Proof.} We drop index $n$ in the course of the proof. Represent $S^i$ as a sum of the process $A^i$ of finite variation and the local martingale $M^i$:
$$ S^i_t=S^i_0+A^i_t+M^i_t,\ \ \ A^i_t=\int_0^t S^i_u\mu^i_u\,du, \ \ M^i_t=\int_0^t (S^i_u\beta^i_u, dW_u).$$
Since $\langle M^i,M^k\rangle_t=\int_0^t S^i_u S^k_u (\beta^i_u,\beta^k_u)du$ we conclude (see \cite{JacShi03}, Theorem 6.30) that the class $L(S)$ consists of predictable processes $\gamma$, verifying the condition
\begin{equation} \label{4.2}
 \int_0^T\left|\sum^{d}_{i=1}\gamma^i_t S^i_t\mu^i_t\right|\,dt+ 
 \int_0^T\sum_{i,k=1}^d \gamma^i_t\gamma^k_t S^i_t S^k_t (\beta^i_t,\beta^k_t)\,dt<\infty\ \ \textnormal{a.s.}
\end{equation}   

Condition (\ref{4.1}) assures that the following local martingale (Girsanov exponential, represented by Dolean-Dade exponential) is well-defined: 
$$ Z_t=\mathcal E(-\lambda\circ W)_t =\exp\left(-\frac{1}{2}\int_0^t |\lambda_u|^2\,du-\int_0^t(\lambda_u,dW_u)\right).$$
By virtue of the representation
$$ S_t^i=S_0^i\,\mathcal E(R^i)_t=S_0^i\exp\left(\int_0^t \left(\mu^i_u-\frac{1}{2} |\beta^i_u|^2 \right) du+\int_0^t(\beta^i_u,dW_u) \right) $$
and the equality $\mu^i=(\beta^i,\lambda)$, we deduce that the processes 
\begin{eqnarray*}
S^i_t Z_t &=& S_0^i \exp\left(-\frac{1}{2}\int_0^t|\beta^i_u-\lambda_u|^2\,du+\int_0^t(\beta^i_u-\lambda_u,dW_u)\right)\\
&=& S_0^i \mathcal E( (\beta^i-\lambda)\circ W)_t
\end{eqnarray*}
are local martingales. It follows (from \cite{ChrLar07}, Lemma 3.2) that 
the processes $XZ$, $X\in\mathcal X$ are local martingales as well.
Hence, they are supermartingales and $Z\in\mathcal D$.

Furthermore, let $\delta=(\sigma\sigma^T)^{-1}\mu$ and $\gamma^i=\delta^i/(ZS^i)$. Then 
$$ \sum_{i=1}^d\gamma^i\mu^i S^i=\frac{1}{Z}(\mu,\delta)=\frac{1}{Z}(\sigma\lambda,\delta)
=\frac{1}{Z}(\lambda,\lambda),$$
$$ \sum_{i,k=1}^d \gamma^i\gamma^k S^i S^k (\beta^i,\beta^k)= \frac{1}{Z^2}\sum_{i,k=1}^d \delta^i\delta^k
(\beta^i,\beta^k)=
 \frac{1}{Z^2}(\sigma^T\delta,\sigma^T\delta)= \frac{1}{Z^2}(\lambda,\lambda)$$
and condition (\ref{4.2}) is satisfied in view of (\ref{4.1}). Thus, $\gamma\in L(S)$. 

Note that
$$(\delta,dR)=(\delta,\mu)dt+(\sigma^T\delta,dW)=(\lambda,\lambda)dt+(\lambda,dW),$$
\begin{eqnarray*}
\mathcal E(\delta\circ R)_t &=& \exp\left((\delta\circ R)_t-\frac{1}{2}\langle \delta\circ R,\delta\circ R\rangle_t\right)\\
&=& \exp\left(\frac{1}{2}\int_0^t |\lambda_u|^2 du+\int_0^t(\lambda_u,dW_u)\right)=\frac{1}{Z_t}.
\end{eqnarray*}
Therefore,  
$$d\left(\frac{1}{Z_t}\right)=\frac{1}{Z_t}(\delta_t,dR_t)=(\gamma,dS)$$
and $Z\in \mathcal D$ admits the representation 
$Z=1/V$, $V=1+\gamma\circ S\in\mathcal X_1$, which was to be proved. \qed

Lemma \ref{lem:2} shows that for the diffusion market model under consideration the results of Sect.~\ref{sec:2} are valid.  The subsequent lemma allows us to give a characterization of NAA and SAA conditions in terms of the sequence $(\lambda^n)$ (see Theorem~\ref{theo:8} below).

\begin{lemma} \label{lem:3}
(a) The following equality holds:
\begin{equation} \label{4.3}
 \lim_{M\to\infty}\limsup_{n\to\infty}\mathsf P^n(V^n_T\ge M)=
 \lim_{M\to\infty}\limsup_{n\to\infty}\mathsf P^n\left(\int\limits_0^{T(n)}|\lambda^n_t|^2 dt\ge M\right);
\end{equation} 
(b) $ \limsup_{n\to\infty}\mathsf P^n(V^n_T\ge M)=1$ for all $M>0$ iff
$$ \limsup_{n\to\infty}\mathsf P^n\left(\int_0^{T(n)}|\lambda^n_t|^2 dt\ge M\right)=1\ \textit{for allõ}\ M>0.$$
\end{lemma}
\textit{Proof.} A simple calculation yields: 
$$\frac{V^n_{T}}{\exp\left(\int_0^{T}|\lambda^n_t|^2 dt\right)} 
=\mathcal E (\lambda^n\circ W^n)_{T},$$
$$ \frac{\exp\left(\frac{1}{2}\alpha (1-\alpha)\int_0^T|\lambda_t^n|^2 dt\right)}{(V_T^n)^\alpha}=
\mathcal E(-(\alpha\lambda^n)\circ W^n)_T,\ \ \alpha\in (0,1).$$
As long as $\mathsf E_{\mathsf P^n}\mathcal (Y^n)_T\le 1$ for any non-negative local martingale $Y^n$, by the estimate  (\ref{1.4}) we obtain:
\begin{equation} \label{4.4}
	\mathsf P^n(V^n_T\ge M)\le \frac{e^N}{M}+\mathsf P^n\left(\int_0^{T}|\lambda^n_t|^2 dt\ge N\right),
\end{equation}
\begin{equation} \label{4.5}
	\mathsf P^n\left(\int_0^{T}|\lambda^n_t|^2 dt\ge M\right) 
	\le\frac{N^\alpha}{e^{\frac{1}{2}\alpha (1-\alpha)M}}+\mathsf P^n(V^n_T\ge N).
\end{equation}

Taking the upper limit as $n\to\infty$, and then the limits as $M\to\infty$, $N\to\infty$, we get equality (\ref{4.3}). Furthermore, let $\limsup_{n\to\infty}\mathsf P^n(V^n_T\ge M)=1$ for all $M>0$. Then inequality (\ref{4.4}) implies that
\begin{equation} \label{4.6}
1\le \frac{e^N}{M} + \limsup_{n\to\infty}\mathsf P^n\left(\int_0^{T}|\lambda^n_t|^2dt\ge N\right).	
\end{equation}
In the case of $\limsup_{n\to\infty}\mathsf P^n\left(\int_0^{T}|\lambda^n_t|^2 dt\ge M\right)=1$ for all $M>0$,
by inequality (\ref{4.5}) we have
\begin{equation} \label{4.7}
	1\le\frac{N^\alpha}{e^{\frac{1}{2}\alpha (1-\alpha)M}}+
\limsup_{n\to\infty}\mathsf P^n(V^n_T\ge N).
\end{equation}
Taking the limits in (\ref{4.6}), (\ref{4.7}) as $M\to\infty$, we get the assertion (b). \qed

\begin{theorem} \label{theo:8}
The following relations hold true:
$$ NAA \iff 
\lim_{M\to\infty}\limsup_{n\to\infty}\mathsf P^n\left(\int_0^{T(n)}|\lambda_t^n|^2 dt
	\ge M\right)=0; $$
$$ SAA \iff \limsup_{n\to\infty}\mathsf P^n\left(\int_0^{T(n)}|\lambda_t^n|^2 dt\ge M\right)=1\ \textit{for all}\ M>0.$$
\end{theorem} 
\textit{Proof} follows directly from Lemma \ref{lem:3} and conditions (c) of Theorems \ref{theo:1} and \ref{theo:2}. \qed

Theorem \ref{theo:8} implies Propositions 8 and 9 of the paper \cite{KabKra98}.
Note, that as compared to \cite{KabKra98}, we do not merely drop the condition  
$\mathsf E_{\mathsf P^n} \mathcal E(-\lambda^n\circ W^n)_T= 1$, but even do not assume the existence of equivalent local martingale measures in the ''small'' markets.

\section{Discrete-time infinite horizon market model with one log-normal stock}
\label{sec:6}
Consider a sequence $(\xi_k)_{k=1}^\infty$ of independent standard normally distributed random variables $\xi_k\in\mathcal N(0,1)$, defined on the probability space $(\Omega,\mathcal F,\mathsf P)$, $\mathcal F=\sigma(\xi_k,\ k\ge 1)$. Assume that there is only one stock, whose price is determined by the recurrence relation 
$$ S_n=S_{n-1}(1+R_n),\ \ \ R_n=\exp(\mu_n-\sigma_n^2/2+\sigma_n\xi_n)-1,\ \ n\ge 1;\ \ S_0=1.$$
Here $\mu_k\in\mathbb R$, $\sigma_k>0$ are non-random sequences. 

We put $\mathbb T^n=\{0,\dots,T(n)\}$, $T(n)=n$ and introduce the sequence of small markets, defined on the probability spaces $(\Omega,\mathcal F_n,\mathsf P)$, $\mathcal F_n=\mathcal F^n_T=\sigma(\xi_1,\dots,\xi_n)$ with associated sets $\mathcal X^n$ of value processes, containing the elements $X$ of the form
$$ X_k=X_{k-1}(1+\delta_k R_k)\ge 0, \ \ k=1,\dots,n;\ \ X_0=x.$$
An element $\delta_k\in [0,1]$ of $(\mathcal F_n)$-predictable process $(\delta_n)$ describes the fraction of wealth, invested in the stock at time $k$.  

Let the elements of the sequence $(\delta_k^*)_{k\ge 1}$ be the solutions of the following optimization problems:
$$ \mathsf E\ln (1+\delta R_k)=
\frac{1}{\sqrt{2\pi}}\int_{-\infty}^\infty\ln\left(1+\delta (e^{\mu_k-\sigma_k^2/2+\sigma_k x}-1)\right)e^{-x^2/2}dx
\to\max_{0\le\delta\le 1}$$
and let $V_n=\prod_{k=1}^n(1+\delta_k^* R_k)$. Then the processes $(X_k/V_k)_{0\le k\le n}$ are supermartingales
(\cite{Bec01}, Example 6). Thus all the results of Sect.~\ref{sec:2} are valid. At that, the measures $\mathsf P^n$ coincide with the restrictions of $\mathsf P$ to $\mathcal F_n$, and $V^n=(V_k)_{k=0}^n$. 

As is shown in \cite{Shi99}, the condition $\sum_{k=1}^\infty(\mu_k/\sigma_k)^2<\infty$ is sufficient for the absence of asymptotic arbitrage. More complete picture is given in Lemma \ref{lem:4} and Theorem~\ref{theo:9}. 
Instead of a direct analysis of the sequence $(V_k)$, we exploit conditions (c) of Theorems \ref{theo:3}, \ref{theo:4}.

Put $\varepsilon>0$ and
$$ \Sigma^1_n(\varepsilon)=\sum_{k=1}^n \left(\frac{\mu_k}{\sigma_k}\right)^2 I_{\left\{0<\mu_k\le \frac{1}{2}(1+\varepsilon)\sigma_k^2\right\}},\ \ \
   \Sigma^2_n(\varepsilon)=\sum_{k=1}^n \mu_k I_{\left\{\mu_k>\frac{1}{2}(1+\varepsilon)\sigma_k^2\right\}},$$ $$\Sigma_n(\varepsilon)=\Sigma^1_n(\varepsilon)+\Sigma^2_n(\varepsilon).$$
\begin{lemma} \label{lem:4}
The following assertions hold:
\begin{itemize}
  \item[(a)] if $\Sigma_\infty(1)<\infty$ then NAA condition is satisfied;
	\item[(b)] if $\Sigma^2_\infty(\varepsilon)=\infty$ for some $\varepsilon>0$ then SAA condition is satisfied;
	\item[(c)] if there exists an $\varepsilon\in (0,1)$ such that	
\begin{equation} \label{5.1}
	\lim_{k\to\infty}\sigma_k I_{\left\{0<\mu_k\le \frac{1}{2}(1+\varepsilon)\sigma_k^2\right\}}=0
\end{equation} 
	and $\Sigma^1_\infty(\varepsilon)=\infty$ then SAA condition is satisfied.
\end{itemize}
\end{lemma}  
\textit{Proof.} (a). Let $Z_n=\prod_{k=1}^n\zeta_k$, where
\begin{eqnarray*}
\zeta_k=I_{\{\mu_k\le 0\}} &+&\exp\left(-\frac{1}{2}\left(\frac{\mu_k}{\sigma_k}\right)^2-\frac{\mu_k}{\sigma_k}\xi_k\right) I_{\{0<\mu_k\le\sigma_k^2\}}\\
&+&\exp\left(-\mu_k+\frac{\sigma_k^2}{2}-\sigma_k\xi_k\right)I_{\{\mu_k> \sigma_k^2\}}.
\end{eqnarray*}

The independence of $\xi_k$ and $\mathcal F_{k-1}$ and the equality $\mathsf E e^{a\xi}=e^{a^2/2}$ imply \begin{eqnarray*}
&&\mathsf E((1 + \delta_k R_k)\zeta_k|\mathcal F_{k-1}) = (1+\delta_k \mathsf E R_k) I_{\{\mu_k\le 0\}}+ \left((1-\delta_k)\mathsf E \zeta_k+\delta_k\right)I_{\{\mu_k>\sigma_k^2\}}\\
&+& \left((1-\delta_k)\mathsf E\zeta_k+	
	\delta_k \mathsf E e^{\left(-\frac{1}{2}\left(\sigma_k-\frac{\mu_k}{\sigma_k}\right)^2+
	\left(\sigma_k-\frac{\mu_k}{\sigma_k}\right)\xi_k\right)}\right)	I_{\{0<\mu_k\le\sigma_k^2\}}\\
&=& \left(1+\delta_k \left (e^{\mu_k}-1\right)\right)I_{\{\mu_k\le 0\}}\\
&+&\left((1-\delta_k)e^{-\mu_k+\sigma_k^2}+\delta_k\right)I_{\{\mu_k>\sigma_k^2\}}
	+ I_{\{0<\mu_k\le\sigma_k^2\}}\le 1.	
\end{eqnarray*}
Thus, $(Z_k)_{k=1}^n$ is an equivalent supermartingale density for $\mathcal X^n$:
$$\mathsf E(X_k Z_k|\mathcal F_{k-1})\le X_{k-1} Z_{k-1},\ \ \ k=1,\dots n;\ \ \ X\in\mathcal X^n.$$

By Theorem \ref{theo:3} we see that NAA condition is a consequence of the equality
$$ \lim_{\alpha\downarrow 0}\liminf_{n\to\infty}\mathsf E Z_n^\alpha=1.	$$
We have,
$$\mathsf E \zeta_k^\alpha =I_{\{\mu_k\le 0\}}+ e^{-\frac{1}{2}\alpha (1-\alpha)\left(\mu_k/\sigma_k\right)^2}
I_{\{0<\mu_k\le\sigma_k^2\}}+e^{-\alpha\mu_k+\frac{1}{2}\alpha (1+\alpha)\sigma_k^2}I_{\{\mu_k>\sigma_k^2\}}.$$
Hence, 
$$\mathsf E Z_n^\alpha\ge \exp\left(-\frac{1}{2}\alpha (1-\alpha)\Sigma^1_\infty(1) -\alpha\Sigma^2_\infty(1)\right)$$
and $\lim_{\alpha\downarrow 0}\liminf_{n\to\infty}\mathsf E Z_n^\alpha\ge 1$ under condition (a). Conversely, $\mathsf E Z_n^\alpha\le 1$, $\alpha\in (0,1)$ by Jensen's inequality.

(b). Put $\delta_k=I_{\left\{\mu_k>\frac{1}{2}(1+\varepsilon)\sigma_k^2\right\}}$ and $\alpha\in (0,\varepsilon)$. Then 
$$\mathsf E(1+\delta_k R_k)^{-\alpha}=e^{-\alpha(\mu_k-\frac{1}{2} (1+\alpha)\sigma_k^2)}\le e^{-\alpha\mu_k \left(1-\frac{1+\alpha}{1+\varepsilon}\right)}\ \	\textnormal{for}\ \ \mu_k>\frac{1}{2}(1+\varepsilon)\sigma_k^2$$
and the correspondent value process $X$ satisfies the inequality
$$ \mathsf E \left(X_n^{-\alpha}\right)=\prod_{k=1}^n \mathsf E(1+\delta_k R_k)^{-\alpha}  \le\exp\left(-\frac{\alpha(\varepsilon-\alpha)}{1+\varepsilon}\Sigma_n^2(\varepsilon)\right). $$
By condition (b) we have $\lim_{n\to\infty}\mathsf E \left( X_n^{-\alpha}\right)=0$ and by Theorem~\ref{theo:4} there exists a strong asymptotic arbitrage.

(c). Put 
$\delta_k=\mu_k\sigma^{-2}_k I_{\left\{0<\mu_k\le \frac{1}{2}(1+\varepsilon)\sigma_k^2\right\}}.$
To show the existence of a strong asymptotic arbitrage it is enough to prove that 
\begin{equation} \label{5.2}
\lim_{n\to\infty}\mathsf E X_n^{-\alpha}=\lim_{n\to\infty}\prod_{k=1}^n \mathsf E(1+\delta_k R_k)^{-\alpha}= 0,\ \ \alpha\in (0,1)
\end{equation}
and apply Theorem~\ref{theo:4}.

Since the series $\Sigma^1(\varepsilon)$ is divergent, the sequence $(\delta_k)$ contains infinitely many non-zero members. Without loss of generality, we assume that $0<\mu_k\le \frac{1}{2}(1+\varepsilon)\sigma_k^2$ for all $k$
and $\sigma_k\to 0$, $k\to\infty$.

By Taylor's formula we obtain the estimate
$$ \psi_\alpha(x)=\frac{1}{(1+x)^\alpha}-(1-\alpha x)\le\frac{\alpha(1+\alpha)}{2}\frac{x^2}{(1+b)^{2+\alpha}},\ \ \
x\ge b>-1.$$ 
Put $\nu_k=-\delta_k\sqrt{\sigma_k}.$ 
Then
\begin{eqnarray*}
\mathsf E \psi_\alpha(\delta_k R_k)&=&\mathsf E\left(\psi_\alpha(\delta_k R_k) 
(I_{\{-\delta_k<\delta_k R_k< \nu_k\}}+I_{\{\delta_k R_k\ge \nu_k\}})\right)\\
&\le &\frac{\alpha(1+\alpha)}{2}\mathsf E\left(\frac{(\delta_k R_k)^2}{(1-\delta_k)^{2+\alpha}}
I_{\{R_k<-\sqrt{\sigma_k}\}}+\frac{(\delta_k R_k)^2}{(1+\nu_k)^{2+\alpha}}I_{\{R_k\ge -\sqrt{\sigma_k}\}}\right)\\
&\le &\frac{\alpha(1+\alpha)}{2}\left(\frac{2^{2+\alpha}\delta_k^2\sigma_k}{(1-\varepsilon)^{2+\alpha}} 
\mathsf P(R_k<-\sqrt{\sigma_k})+\frac{\delta_k^2\mathsf E(R_k^2)}{(1+\nu_k)^{2+\alpha}}\right).
\end{eqnarray*}
Let us find an asymptotic form of the right-hand side of this inequality as $k\to\infty$. By formula (7.1.23) of \cite{AbrSte72}:
\begin{eqnarray*}
\mathsf P(\xi_k<-x)
&=&\frac{1}{\sqrt{2\pi}}\int_{-\infty}^{-x} e^{-t^2/2}\,dt
=\frac{1}{\sqrt{\pi}}\int_{x/\sqrt{2}}^{+\infty} e^{-\tau^2}\,d\tau\\
&=&\frac{1}{2}\erfc\left(\frac{x}{\sqrt 2}\right)\sim \frac{e^{-x^2/2}}{\sqrt{2\pi}x},\ \ x\to+\infty,
\end{eqnarray*}
we get
\begin{eqnarray*}
\mathsf P(R_k<-\sqrt{\sigma_k})=\mathsf P\left(\xi_k< \frac{\sigma_k}{2}-\frac{\mu_k}{\sigma_k}+\frac{1}{\sigma_k}\ln\left(1-\sqrt{\sigma_k}\right)\right)
\\ \le \mathsf P\left(\xi_k<\frac{\sigma_k}{2}-\frac{1}{\sqrt{\sigma_k}}\right) 
\sim 
\frac{\sqrt{\sigma_k}}{\sqrt{2\pi}}\exp\left(-\frac{1}{2\sigma_k}\right),\ \ k\to\infty.
\end{eqnarray*}
Moreover, $\nu_k\to 0$ and
$$\mathsf E R_k^2=e^{2\mu_k+\sigma_k^2}-2e^{\mu_k}+1\sim\sigma_k^2,
 \ \ k\to\infty.$$
Consequently,
$$ \mathsf E\psi_\alpha(\delta_k R_k)\le\frac{\alpha(1+\alpha)}{2}A_k,\ \ \ 
A_k\sim\delta_k^2\sigma_k^2=\left(\frac{\mu_k}{\sigma_k}\right)^2,\ \ k\to\infty.$$

Furthermore, $\mathsf E R_k=e^{\mu_k}-1\ge\mu_k$. Thus,
\begin{eqnarray*}
\mathsf E\frac{1}{(1+\delta_k R_k)^\alpha}=1-\alpha\delta_k\mathsf E R_k+\mathsf E\psi_\alpha(\delta_k R_k)
\le 1-\alpha\left(\frac{\mu_k}{\sigma_k}\right)^2+\frac{\alpha(1+\alpha)}{2}A_k\\
\le\exp\left(-\alpha\left(\frac{\mu_k}{\sigma_k}\right)^2\left(1-\frac{(1+\alpha)}{2}(1+o(1))\right)\right),\ \
 k\to\infty
\end{eqnarray*}
and the condition $\Sigma_\infty^1(\varepsilon)=\infty$ implies (\ref{5.2}). \qed

\begin{theorem} \label{theo:9}
The following assertions hold:
\begin{itemize}
\item[(a)] exactly one of the conditions NAA or SAA is satisfied and 
  $$ NAA\iff\mathsf P(V_\infty<\infty)=1,\ \ \ SAA\iff\mathsf P(V_\infty=\infty)=1.$$
\item[(b)] If condition (\ref{5.1}) is satisfied for some $\varepsilon\in (0,1)$ then
$$ NAA \iff \Sigma_\infty(\varepsilon)<\infty;\ \ \ SAA \iff \Sigma_\infty(\varepsilon)=\infty.$$	 
\end{itemize}
\end{theorem}
\textit{Proof.} The process $(1/V_n)_{n\ge 1}$ converges a.s. since it is a positive supermartingale. Moreover, $V_n$ is a product of independent positive random variables. From Kol\-mo\-go\-rov's zero-one law it follows that the events $\{V_\infty<\infty\}$, $\{V_\infty=\infty\}$ have probability $0$ or $1$. 

Let $\mathsf P(V_\infty<\infty)=1$. The sequence $V_n$ converges to $V_\infty<\infty$ a.s. Hence, $(V_n|\mathsf P)$ is tight:
$$\limsup_{n\to\infty} \mathsf P(V_n\ge M)\le\mathsf E\left( \limsup_{n\to\infty} I_{\{V_n\ge M\}}\right)
\le\mathsf P(V_\infty\ge M/2),\ \ M>0.$$
By Theorem~\ref{theo:1} this implies the absence of asymptotic arbitrage.

If $\mathsf P(V_\infty=\infty)=1$ then 
$$\limsup_{n\to\infty} \mathsf P(V_n< M)\le\mathsf E\left( \limsup_{n\to\infty} I_{\{V_n< M\}}\right)
=0,\ \ M>0$$
and condition (c) of Theorem~\ref{theo:2} is satisfied. Hence, there exists a strong asymptotic arbitrage.

Assertion (b) follows from Lemma \ref{lem:4} and the inequality $\Sigma_n(\varepsilon)\ge\Sigma_n(1)$, $\varepsilon\in (0,1)$. \qed




\end{document}